\documentclass{amsart}
\usepackage{amscd,amssymb,euscript}

\newcounter{noindnum}[subsection]
\newcommand{\noindstep}{\refstepcounter{noindnum}{\rm(}\arabic{noindnum}\/{\rm) }}

\renewcommand{\phi}{\varphi}
\renewcommand{\epsilon}{\varepsilon}

\newcommand{\divisor}{\mathfrak D}
\newcommand{\Dmod}{\EuScript D}

\newcommand{\D}{\mathcal D^\bullet}
\newcommand{\K}{\mathcal K^\bullet}

\newcommand{\Oo}{\mathcal O}
\newcommand{\T}{\mathcal T}

\newcommand{\U}{\mathcal U}
\newcommand{\V}{\mathcal{V}}
\newcommand{\W}{\mathcal{W}}

\newcommand{\C}{\mathbb C}
\newcommand{\Pro}{\mathbb P}

\newcommand{\Z}{\mathbb{Z}}
\newcommand{\Hyper}{\mathbb H}

\newcommand{\M}{\mathbf M}
\newcommand{\A}{\mathbf A}
\newcommand{\X}{\mathbf X}
\newcommand{\Y}{\mathbf Y}

\newcommand{\Id}{\mathop{\mathrm{Id}}\nolimits}
\newcommand{\Spec}{\mathop{\mathrm{Spec}}}
\newcommand{\Conn}{\lefteqn{\widetilde{\phantom{Conn}}}}

\theoremstyle{plain}
\newtheorem*{theorem*}{Theorem}
\newtheorem{theorem}{Theorem}
\newtheorem{proposition}{Proposition}
\newtheorem{lemma}{Lemma}

\theoremstyle{definition}
\newtheorem{definition}{Definition}

\theoremstyle{remark}
\newtheorem*{remark}{Remark}
\newtheorem*{remarks}{Remarks}

\title{Frobenius manifold structures\\
on the spaces of abelian integrals}
\author{Roman M.\ Fedorov}
\email{fedorov@bu.edu}
\date{}
\address{Department of Mathematics and Statistics\\
Boston University\\
111 Cummington St\\
Boston, MA, 02215}

\begin{document}

\begin{abstract}
Frobenius manifold structures on the spaces of abelian integrals were
constructed by I.~Krichever. We use $\Dmod$-modules, deformation theory, and
homological algebra to give a coordinate-free description of these
structures. It turns out that the tangent sheaf multiplication has a
cohomological origin, while the Levi--Civita connection is related to
one-dimensional isomonodromic deformations.
\end{abstract}

\maketitle

\thispagestyle{empty}

\section{Introduction}
Frobenius manifolds are manifolds with a flat metric and a multiplication in
the tangent sheaf, subject to some constraints. Frobenius manifolds were
introduced by B.~Dubrovin in~\cite{Dubrovin1,Dubrovin2} as a mathematical framework for
deformations of topological quantum field theories (see also~\cite{Dubrovin}). In mathematics
Frobenius manifolds arise in two different situations, corresponding to
A-models and B-models in physics. In an A-model one counts rational curves on
a variety; this is also known as Gromov--Witten invariants. The generating
function for these invariants is the potential for the corresponding
Frobenius manifold.

This paper is concerned with B-models. In a B-model one studies deformations
of a certain complex structure (formal or analytic). The best known examples
are extended moduli spaces of Calabi--Yau varieties~\cite{KonBar} and the
unfoldings of isolated singularities~\cite{Saito} (see~\cite{Sabbah} for an
exposition). We would like to mention that Frobenius structures are important
for mirror symmetry: if two varieties are mirror dual to each other, then the
A-model Frobenius manifold, corresponding to the first variety, is isomorphic
to the B-model Frobenius manifold, corresponding to the second.

\subsection{Moduli spaces of abelian integrals}
Examples of Frobenius manifolds are furnished by Hurwitz spaces.
Hurwitz spaces parameterize pairs $(X,f)$, where $X$ is a smooth
complete algebraic curve, $f:X\to\Pro^1$. Dubrovin constructed
Frobenius structures on Hurwitz spaces~\cite{Dubrovin}.

Our main object is the following deformation of a Hurwitz space: a space of
pairs $(X,f)$, where $f$ is a \emph{multi-valued} function such
that $df$ is a single-valued meromorphic 1-form with prescribed periods and residues. If
the periods and residues are equal to zero, then this space is a
Hurwitz space. Our spaces will be called \emph{spaces of abelian integrals}.

Krichever constructs in~\cite{Krichever1,Krichever2} Frobenius structures on
the universal covers of the spaces of abelian integrals. Our main goal is to
give a coordinate-free geometric description of these Frobenius structures. We also generalize the setup to the case of multiple poles and non-zero residues in~\S\ref{sectmulti}. In particular, our generalization covers the previously untreated case of abelian integrals of the third kind.
Our approach is based on a $\Dmod$-module push-forward (also known as twisted
de Rham complex; see~\cite[\S I.3.3]{Sabbah}). It turns out that these
structures of Frobenius manifolds have a nice interpretation: the tangent
sheaf multiplication has a cohomological origin, similar to that
of~\cite{KonBar}. The metric and the Levi--Civita connection are closely
related to one-dimensional isomonodromic deformation. (This is not directly
related to isomonodromic deformations used to describe the semi-simple
Frobenius manifolds.)

We are using the approach to Frobenius structures via primitive forms. This
has been invented by K.~Saito~\cite{Saito}.  We would like to mention a
striking similarity between three constructions of Frobenius structures:
on the universal unfolding of isolated singularity~\cite{Saito}, on the
extended moduli space of Calabi--Yau varieties~\cite{KonBar}, and our
construction. In each case a pencil of connections is obtained by the
(derived) direct image. Our case is, in some sense, intermediate: on the one
hand, singularities are present, on the other hand, our structure is not
local, it depends on the global geometry of a curve. This is why we hope
that, generalizing our construction to higher dimension, we shall provide a
bridge between the pictures of Saito and Barannikov--Kontsevich, giving a
unified approach to B-model Frobenius manifolds.

Another interesting feature of our construction is that we get a
\emph{family\/} of Frobenius manifolds parameterized by the periods of
abelian integrals. We also want to emphasize that there are some new features
specific for the higher genus case: to get a Frobenius structure we need to
make a modification of the direct image (see~\S\ref{improving}).

My interest in the subject was attracted by
papers~\cite{Krichever2,Losev}. I would like to thank
D.~Arinkin, V.~Ginzburg, A.~Losev, and I.~Mirkovi\'c for
numerous discussions. A part of this work was done, while I
was visiting Max Planck Institute in Bonn. I would like to thank it
for warm hospitality. I am also grateful to the reviewer for the valuable comments.

\section{Preliminaries and the main construction}
\subsection{Pencils of connections}

\begin{definition}\label{pencil}
Let $p_1:\M\times\Pro^1\to\M$ be the natural projection, where $\M$ is a
complex manifold, $\Pro^1=\Pro^1_\C$. By a \emph{pencil of connections\/} on $\M$ we mean a pair
$(\W,\nabla)$, where $\W\to\M$ is a vector bundle, $\nabla$ is a relative
flat connection on $\V=p_1^*\W$ along $\M$ with a simple pole along
$\M\times\{0\}$.
\end{definition}

One interprets a pencil as a family of flat connections on $\W$,
parameterized by $\Pro^1\setminus0$. The condition on the pole implies that
this family is of the form $\nabla_\infty+\Phi/z$, this is why it is called
``pencil'' (here $z$ is a coordinate on $\Pro^1$). There is a natural way to construct twisted Frobenius manifold
structures on dense open subsets of $\M$ starting from a pencil
of connections, provided this pencil of connections satisfies
some non-degeneracy condition (and, conversely, every Frobenius
manifold gives rise to a pencil of connections). This will be
explained in details in~\S\ref{WDVV}.

\subsection{Main objects}\label{mainobj}
Consider a smooth complete algebraic curve $X$ of genus $g$ over $\C$, let
$p\in X$. Denote by $(\hat X,\hat p)$ the maximal abelian cover of $(X,p)$.

\begin{definition}
An \emph{abelian integral\/} on $(X,p)$ is a function $f$ on
$\hat X$ such that $df$ descends to a meromorphic differential on $X$. We
define \emph{periods\/} of $f$ to be those of $df$.
\end{definition}

\begin{remarks}
(1) An abelian integral can be thought as an integral of a meromorphic form
on $X$. Thus the space of abelian integrals is a one-dimensional affine bundle
over the vector bundle of differential forms.

(2) One can avoid working with abelian integrals by fixing a point
$p_0\in\hat X$. Then the affine bundle trivializes, a section being the space
of abelian integrals that vanish at $p_0$, see Remark in~\S\ref{sectmulti}.
\end{remarks}

To simplify notation, we shall assume first that $df$ has a single pole. We outline the changes needed in the multi-pole case in~\S\ref{sectmulti}.

Let $n\ge1$ be an integer. Consider the moduli space $\A_{g,n}$ of triples $(X,p,f)$, where $(X,p)$ is
as above, $f$ is an abelian integral with a single pole of order exactly $n$ at $p$
(in other words, $df$ is a meromorphic form on $X$ with the only pole at $p$
of order $n+1$).

The periods of $f$ give a linear map $H_1(X,\Z)\to\C$. One can identify
groups $H_1(X,\Z)$ locally over the moduli space of curves using the
Gauss--Manin connection, therefore the periods give rise to a foliation on
$\A_{g,n}$. Let us fix one of the leaves and denote its smooth locus by $\A$. Thus, roughly speaking, $\A$
parameterizes abelian integrals with prescribed periods.

Let $\hat\A$ be the moduli space of quadruples $(X,p,f,\Delta)$, where
$(X,p,f)\in\A$, $\Delta$ is a subgroup of $H_1(X,\Z)$ maximal isotropic with
respect to the intersection form. The elements of $\Delta$ will be called
\emph{$a$-cycles}. Clearly, $\hat\A$ is a cover of $\A$. Our main result is
the following
\begin{theorem*}
{\rm(}a{\rm)} There is a natural pencil of connections on $\hat\A$.\\
{\rm(}b{\rm)} This pencil of connections gives rise to a twisted Frobenius
structure on $\hat\A$.
\end{theorem*}

Let us comment on the second part. To construct a Frobenius manifold
structure from the pencil of connections one needs a so-called
\emph{primitive section\/}, which exist only locally. Another problem is that
they are not canonical (see~\S\ref{WDVV}). It means that~$\hat\A$ is covered
by ``Frobenius charts''. The tangent sheaf multiplication is nevertheless
canonical. This is what we mean by a twisted Frobenius structure. It will be
discussed in more details in~\S\ref{WDVV}.
We shall further discuss the choice of a primitive section in~\S\ref{PrimitiveSections}.

\begin{remarks}
(1) We could have worked with fibrations rather than foliations by fixing
a~basis in $H_1(X,\Z)$ (that is, a level structure on the curve). Then the
periods give a global map to $\C^{2g}$ and our moduli spaces are the fibers
of this map.

(2) $\A$ is not algebraic, so we shall work in analytic category. Also, we can construct a Frobenius structure on
the formal completion of $\A$ at a point. However, we wanted to emphasize
that our construction is global rather than formal.
\end{remarks}

\subsection{Notation}\label{notation}
The following notation will be fixed throughout the paper. Let $\phi:\X\to\A$
be the universal curve, that is, the fiber of $\phi$ over $(X,p,f)$ is $X$. We
denote by $d_X$ the relative differential $\Oo_\X\to\Omega_{\X/\A}$. We
denote by $d:\Oo_\X\to\Omega_\X$ the usual (absolute) differential.

There is a natural section of $\phi$ corresponding to $p\in X$, denote it by
$\tilde p$\,; we can also view it as a divisor on $\X$. For any integer
number $k$ set $\Oo(k)=\Oo_\X(k\tilde p)$,
$\Omega(k)=\Omega_{\X/\A}\otimes\Oo(k)$, $\T(k)=\T_{\X/\A}\otimes\Oo(k)$, where $\Omega_{\X/\A}$ is the sheaf of relative differentials, $\T_{\X/\A}$ is the relative tangent sheaf.

We have a ``universal'' multi-valued function on $\X$. We denote it again by
$f$, hopefully it will not lead to a confusion. Set $\omega=d_Xf$.

\subsection{Main construction}\label{MainConstruction}
The fact that the periods of $f$ are fixed shows that $df$ is a single-valued
1-form on $\X$. Thus
\begin{equation*}
       \nabla=d+\frac{df}z
\end{equation*}
is a family of flat connections on $\Oo_\X$ parameterized by
$z\in\Pro^1\setminus0$. We can view~$\nabla$ as a relative $\Dmod$-module on
$\X\times(\Pro^1\setminus0)$. The idea is that we get a pencil of connections
on $\A$ by taking the push-forward of $\nabla$ along
$\phi\times\Id_{\Pro^1}$. There are two problems we shall have to go around

\noindstep Our $\Dmod$-module is not defined at $z=0$, and the push-forward
is not coherent near $z=\infty$. Thus some regularization is needed.

\noindstep
We get a vector bundle on $\A\times\Pro^1$, whose
restriction to
$\{m\}\times\Pro^1$ is not trivial, thus it is not a pencil of connections.
We shall make some modification along $\A\times\{\infty\}$, this
is where we need the additional structure of $a$-cycles.

We shall denote $\phi\times\Id_{\Pro^1}$,
$\phi\times\mathrm\Id_{\Pro^1\setminus0}$ etc. again by $\phi$ for brevity.
\setcounter{noindnum}{0}
\begin{remarks}
\noindstep\label{FirstRemark} For a fixed $z$ we can start with a connection
$d_X+\frac\omega{z}$ along the fibers of $\phi$. Then the
condition that the periods of $\omega$ are constant is exactly
the isomonodromic condition for this connection. Thus it can be
extended to an absolute connection.
See~\S\ref{isomonodromy} for more on isomonodromic deformation.

\noindstep If $f$ has zero periods (that is, $f$ is a meromorphic function on
$X$), then one can extend $\nabla$ to an absolute flat meromorphic connection
on $\X\times\Pro^1$. In this case we get a Frobenius manifold with Euler
field (that is, with a conformal structure). For details see, e.g.~\cite{Sabbah}.

\noindstep The $\Dmod$-module push-forward can be viewed as taking the
cohomology fiberwise, with the connection on cohomology being the
Gauss--Manin connection. One can also think about this as a de Rham
complex, twisted by $e^f$; see, e.g.~\cite[\S I.3.3]{Sabbah}.
\end{remarks}

\subsection{Organization of the paper}
In the next section we shall make the ideas above precise, thus constructing
a pencil of connections on $\hat\A$. In~\S\ref{IdenStr} we prove that the
pencil of connections gives rise to Frobenius structures on some open subsets
of $\hat\A$ and calculate these Frobenius structures explicitly, then we
present the relation between our construction and that
of~\cite{Krichever1,Krichever2}. In~\S\ref{sectmulti} we
generalize our setup to the multi-pole case. In~\S\ref{PrimitiveSections}
we discuss the choices involved in the construction of Frobenius manifolds from the pencil of connections.

\section{The precise construction}
Consider the complex
\begin{equation*}
    \Oo(-1)\xrightarrow{d_X+\frac\omega{z}}\Omega(n).
\end{equation*}
We always place the leftmost term in degree zero. We can view this complex as
a complex of sheaves on $\X$, depending on a parameter
$z\in\Pro^1\setminus0$. We can also view it as a complex of
sheaves of $\X\times(\Pro^1\setminus0)$-modules
$\Oo(-1)\boxtimes\Oo_{\Pro^1\setminus0}\xrightarrow{d_X+\omega/z}
\Omega(n)\boxtimes\Oo_{\Pro^1\setminus0}$.

To regularize the complex at $z=0$ we shall patch it on
$\Pro^1\setminus\{0,\infty\}$ with the complex
$\Oo(-1)\xrightarrow{zd_X+\omega}\Omega(n)$, using the diagram
\begin{equation}\label{patch}
\begin{CD}
\Oo(-1)@>zd_X+\omega>>\Omega(n)\\
@V=VV @V\times\frac1zVV\\
\Oo(-1)@>d_X+\frac\omega z>>\Omega(n)\\
\end{CD}
\end{equation}
This is an isomorphism of complexes if $z\ne0,\infty$. This construction
gives a complex
\begin{equation}\label{global}
    \Oo(-1)\boxtimes\Oo_{\Pro^1}\xrightarrow{d_X+\frac\omega z}
    \Omega(n)\boxtimes\Oo_{\Pro^1}(1).
\end{equation}
We shall denote this complex by $\D$ and its restriction to $\X\times\{z\}$
by $\D_z$.

Note that the map is $\Oo_{\A\times\Pro^1}$-linear. We shall view this
complex as an object in the derived category of complexes of sheaves of
$\phi^\cdot\Oo_{\A\times\Pro^1}$-modules on $\X\times\Pro^1$. Here
$\phi^\cdot$ is the sheaf-theoretic inverse image. We shall abuse the
language by saying ``complex of sheaves'', where we really mean ``complex of
sheaves up to a quasi-isomorphism''. The push-forward $\phi_*\D$ is a complex
of coherent sheaves (because $\phi$ is proper).

\subsection{Relation to $\Dmod$-modules}\label{Dmodules}
Denote by $\Oo(\infty)$ the sheaf of meromorphic functions on $\X$ with poles on $\tilde p$ of any order.
The connection $d+\frac{df}z$ equips $\Oo(\infty)$ with a structure of
$z$-dependent $\Dmod_{\X}$-module ($z\ne0$).

To calculate its $\phi$-push-forward we have to consider the corresponding
relative de Rham complex
\begin{equation}\label{Dmod}
    \Oo(\infty)\xrightarrow{d_X+\frac\omega z}\Omega(\infty).
\end{equation}
It will be more convenient for us to work with coherent sheaves, so we take a
subcomplex of coherent sheaves
\begin{equation*}
    \Oo(k)\xrightarrow{d_X+\frac\omega z}\Omega(k+n+1).
\end{equation*}
One easily checks that it is quasi-isomorphic to~(\ref{Dmod}) at
$z\ne\infty$ for any $k\in\Z$. It follows that $\phi_*\D$ has a $\Dmod_\A$-module structure
for $z\ne0,\infty$ (because it is isomorphic to a push-forward of a
$\Dmod$-module). We shall calculate this structure and see that it extends
to $z=\infty$ and has a simple pole at $z=0$.

Our choice of $k=-1$ is imposed by the fact that in this case the
push-forward is locally free, as we shall see shortly. If, for example, we
take $k=0$, the push-forward would have torsion at $z=\infty$ and would not
be concentrated in a single homological dimension.

\subsection{Study of $\phi_*\D$ in the direction of $\Pro^1$}\label{32}
Now we choose a point $m=(X,p,f_0)\in\A$ and study the
restriction of $\phi_*\D$ to $\{m\}\times\Pro^1$. We denote this restriction
again by $\phi_*\D$ for brevity. In this and in the next subsections $\Oo(k)$ stands for
$\Oo_X(kp)$ and $\Omega(k)=\Omega_X(kp)$.

\begin{lemma}\label{structure}
$\phi_*\D$ is isomorphic to the vector bundle
\begin{equation}\label{isom}
\left(\oplus_1^g\Oo_{\Pro^1}\right)\oplus
\left(\oplus_1^{g+n-1}\Oo_{\Pro^1}(1)\right)
\end{equation}
placed in degree~1.
\end{lemma}
\begin{proof}
Let us show first that $\phi_*\D$ is a vector bundle concentrated in
degree~1. To this end we fix $z\ne0$ and consider the exact sequence
\begin{equation}\label{exact}
    0\to\Omega(n)[1]\to\D_z\to\Oo(-1)\to0.
\end{equation}
Let us write the corresponding exact sequence of hypercohomology
\begin{equation*}
\small{0\to\Hyper^0(\D_z)\to0\to
    H^0(\Omega(n))\to\Hyper^1(\D_z)\to H^1(\Oo(-1))\to0\to\Hyper^2(\D_z)\to0.}
\end{equation*}
We see that the hypercohomology of $\D_z$ are concentrated in degree~1 and
the dimension does not depend on $z$ (it is equal to $2g+n-1$). A similar
sequence for the upper complex in~(\ref{patch}) shows that a similar
statement is valid near $z=0$. Now a standard argument involving a base change shows
that $\phi_*\D$ is a locally
free sheaf in degree~1, that is, a vector bundle.

Evaluating the global sections of~(\ref{global}) along $\Pro^1$ first, we
come to the following presentation of the global sections of $\phi_*\D$
\begin{equation*}
    \Hyper^1(\Oo(-1)\xrightarrow{(d_X,\times\omega)}\Omega(n)\oplus\Omega(n)).
\end{equation*}
Using an exact sequence, similar to~(\ref{exact}), one checks easily that the
dimension of the space of global sections of $\phi_*\D$ is equal to
$3g+2n-2=\dim\phi_*\D+(g+n-1)$. The lemma will be proved if we show that (i)
the global sections generate the fiber of $\phi_*\D$ at $z=0$ and (ii) that
there are no global sections, vanishing at both $z=0$ and $z=\infty$. Indeed,
every vector bundle on $\Pro^1$ is isomorphic to $\oplus_i\Oo_{\Pro^1}(m_i)$.
Now (i) shows that $m_i\ge0$ for all $i$, and (ii) shows that $m_i\le1$ for
all $i$.

It follows from the base change that the following map of complexes gives
rise to the map from the space of global sections to $\phi_*\D|_{z=0}$
\begin{equation}\label{restriction}
\begin{CD}
    \Oo(-1)@>(d_X,\times\omega)>>\Omega(n)\oplus\Omega(n)\\
    @VVV @V(0,1)VV\\
    \Oo(-1)@>\omega>>\Omega(n).
\end{CD}
\end{equation}
The second hypercohomology group of the kernel of this map is isomorphic to
$H^1(\Omega(n))=0$. Hence the induced map on hypercohomology groups is
surjective, and (i) is satisfied.

We also see that the kernel of the map from the space of global sections to
$\phi_*\D|_{z=0}$ is given by the global sections of the \emph{first\/}
$\Omega(n)$ summand.

Now, one writes a map of complexes, analogous to~(\ref{restriction}) but
generating a map from the space of global sections to $\phi_*\D|_{z=\infty}$
and checks that its kernel is given by the global sections of the
\emph{second\/} $\Omega(n)$ summand. Since these subspaces of global sections
do not intersect, (ii) follows.
\end{proof}
Below we often view $\phi_*\D$ as a vector bundle, making the shift by $-1$ implicit.

Let us denote by $Conn^0_{X,p,n}$ the moduli space of degree zero line bundles on
$X$ with a connection such that the connection has a pole of order at most
$n$ at $p$ and no other poles. Denote by $\Conn Conn^0_{X,p,n}$ the universal
cover of $Conn^0_{X,p,n}$. Note that $Conn^0_{X,p,n}$ is an
affine bundle over $Pic^0_X$.

\begin{lemma}\label{isomorhism}
There is a natural isomorphism
\begin{equation}\label{identif}
    \phi_*\D|_{z=\infty}\approx\Conn Conn^0_{X,p,n}.
\end{equation}
\end{lemma}
\begin{proof}
By the base change LHS is given by $\Hyper^1(\Oo(-1)\xrightarrow
d\Omega(n))$. It is easy to see that the natural inclusion of the above
complex into $\Oo\xrightarrow d\Omega(n)$ induces an isomorphism in the first
hypercohomology groups. Further, $\Hyper^1(\Oo\xrightarrow d\Omega(n))$ is
identified with the tangent space to $Conn^0_{X,p,n}$ at zero. The latter
space is identified with the universal cover of $Conn^0_{X,p,n}$.
\end{proof}
Below we always assume the identification~(\ref{identif}).

\subsection{Improving $\phi_*\D$}\label{improving}
According to Definition~\ref{pencil}, in order to give rise to a pencil of
connections, $\phi_*\D$ has
to be isomorphic to $p_1^*\W$, where $\W$ is a vector bundle on $\A$. This is
impossible, since the restriction of $\phi_*\D$ to $\{m\}\times\Pro^1$ is not
a trivial bundle (see Lemma~\ref{structure}). This is easy to cure if $g=0$:
just twist by $\Oo_{\Pro^1}(-1)$. If $g>0$ we need to choose a trivial
subbundle in the restriction of $\phi_*\D$ to every $m\in\A$.

Let us now recall that we have chosen a space $\Delta$ of $a$-cycles on $X$.
Note that every degree zero line bundle has a unique non-singular connection
with trivial $a$-monodromy. This gives a splitting of the standard exact
sequence
\begin{equation}\label{UsualExact}
    0\to H^0(X,\Omega(n))\to Conn^0_{X,p,n}\to Pic^0_X\to0.
\end{equation}
Thus, the choice of $\Delta$ gives a splitting
\begin{equation}\label{splitting}
    \phi_*\D|_{z=\infty}=H^0(X,\Omega(n))\oplus H^1(X,\Oo).
\end{equation}
We define $\V$ as the subsheaf of $\phi_*\D$ whose sections over
$U\subset\Pro^1$ are given by
\[
\Gamma(U,\V)=\{s\in\Gamma(U,\phi_*\D):\;s|_{z=\infty}\in H^1(X,\Oo)\}.
\]
\begin{proposition}
The bundle $\V$ is trivial.
\end{proposition}
\begin{proof}
Note that the isomorphism~\eqref{isom} is not canonical but the inclusion of $\V':=\oplus_1^{g+n-1}\Oo_{\Pro^1}(1)$ into $\phi_*\D$ is canonical. It is easy to see that the proposition is equivalent to $\V'_\infty\cap H^1(X,\Oo)=0$. Thus it is enough to prove that $\V'$ corresponds to the connections on the trivial bundle under~\eqref{identif}. To prove the claim, we note that
\[
    \V'_\infty=\{s(\infty)|s\in H^0(\Pro^1,\phi_*\D), s(0)=0\}.
\]
Now the claim follows from the proof of Lemma~\ref{structure} (see diagram~\eqref{restriction} and two paragraphs after it).
\end{proof}

\subsection{Globalization}
So far we were working with a fixed quadruple $(X,p,f_0,\Delta)$. Let us now
globalize the picture to $\hat\A$. First we get a vector bundle $\phi_*\D$ on
$\A\times\Pro^1$ with a connection along $\A$ with a pole along
$\A\times\{0\}$. The restriction of this connection to $\A\times\{z\}$ we
denote by $\nabla_z$. Modifying $\phi_*\D$ at $z=\infty$ as above we get a
vector bundle $\V$ on $\hat\A\times\Pro^1$, its restriction to any point of
$\hat\A$ being a trivial bundle on $\Pro^1$. Thus $\V=p_1^*\W$ for some
vector bundle $\W$ on $\hat\A$.
\begin{theorem}\label{identifications}
$(\W,\nabla)$ is a pencil of connections on $\hat\A$ in the sense of
Definition~\ref{pencil}.
\end{theorem}
\begin{proof}
As was explained in~\S\ref{Dmodules}, $\phi_*\D$ is equipped with a
connection $\nabla_z$ for $z\ne0,\infty$, thus $\V$ also has a $z$-dependent
connection. It is not hard to show directly that this connection extends to
$z=\infty$ and has a simple pole at $z=0$. However, in the next section we
shall calculate both the residue of $\nabla_z$ at zero and $\lim_{z\to\infty}\nabla_z$
explicitly, completing the proof of the theorem.
\end{proof}

\section{Identifying structures}\label{IdenStr}
We shall see shortly that Frobenius structures come from $\nabla_\infty$ and
the residue of the pencil of connections. Then we shall calculate these parts
of the pencil explicitly. It will yield Theorem~\ref{identifications}. It
will also follow that this pencil of connections gives rise to twisted Frobenius
manifold structures on $\hat\A$.

\subsection{From pencils of connections to Frobenius manifolds and WDVV
equation}\label{WDVV}
Let $(\W,\nabla)$ be a pencil of connections on a manifold $\M$. One
interprets it as a family of flat connections on $\M$, parameterized by
$\Pro^1\setminus0$. Our condition on the pole at zero implies that this
family is of the form $\nabla_\infty+\Phi/z$, where $z$ is the standard
coordinate on $\Pro^1$, $\nabla_\infty$ is the restriction of $\nabla$ to
$\M\times\{\infty\}$, $\Phi$ is a Higgs field, that is, an
\emph{$\Oo_\M$-linear} map $\W\to\W\otimes\Omega_\M$.

Assume now that there exists a \emph{primitive section\/}, i.e., a section
$\rho$ of $\W$ such that $\nabla_\infty\rho=0$ and
$\Phi(\rho,\cdot):\T\M\to\W$ is an isomorphism. One uses this isomorphism to
carry $\nabla_\infty$ and $\Phi$ to $\T\M$. The former becomes a flat
structure $\tilde\nabla$ on $\M$ (it is automatically without torsion). The
latter becomes a commutative associative multiplication in the tangent sheaf.
We denote this multiplication by $\circ$. The equation $\Phi(\rho,e)=\rho$
defines a unit for $\circ$. One can show that $\circ$ does not depend on the
choice of $\rho$ (while $\tilde\nabla$ does depend).

The last ingredient needed to equip $\M$ with a structure of a Frobenius
manifold is a symmetric bilinear product compatible with $\tilde\nabla$ and $\circ$. In
other words, $\tilde\nabla$ is the Levi--Civita connection for this metric,
while the multiplication operators are symmetric. These structures altogether
make $\T\M$ into a sheaf of \emph{Frobenius algebras}, so that~$\M$ becomes a
\emph{Frobenius manifold}.

Set $a(z)=-z$, then giving such a metric is equivalent to a $\nabla$-flat
$a$-symmetric non-degenerate pairing $\langle\cdot,\cdot\rangle:\V\otimes
a^*\V\to\Oo_{\M\times\Pro^1}$, where $\V$ is the pull-back of $\W$ to $\M\times\Pro^1$. Indeed, if $s_1$ and $s_2$ are sections
of~$\W$, then we have
\begin{equation*}
    \left\langle\left(\nabla_\infty+\frac\Phi z\right)s_1,s_2\right\rangle+
\left\langle s_1,\left(\nabla_\infty-\frac\Phi
z\right)s_2\right\rangle=d\left\langle s_1,s_2\right\rangle.
\end{equation*}
This is equivalent to $\langle\nabla_\infty s_1,s_2\rangle+ \langle
s_1,\nabla_\infty s_2\rangle=d\langle s_1,s_2\rangle$ and $ \langle\Phi
s_1,s_2\rangle=\langle s_1,\Phi s_2\rangle$.

One then chooses flat coordinates $t_A$ on $\M$ and shows that there exists
locally on $\M$ a \emph{potential} $F$ such that
\begin{equation*}
    \langle\partial_A\circ\partial_B,\partial_C\rangle=
    \tilde\nabla_{\partial_A}\tilde\nabla_{\partial_B}\tilde\nabla_{\partial_C}F.
\end{equation*}
The associativity condition for $\circ$ transforms into WDVV equation for
$F$.

\begin{remark}
In the original definition of Frobenius manifold the structures above are
required to be homogeneous, so the Euler field is added to the structures.
Thus our definition is more general, in fact we shall construct non-homogeneous
Frobenius manifolds.
\end{remark}

\subsubsection{Constructing primitive sections}\label{constructing} Let $m$ be a point of $\M$,
$\rho_m\in\W_m$ be such that $\Phi(\rho_m,\cdot):\T_m\M\to\W_m$ is an
isomorphism. Then we can extend $\rho_m$ to a $\nabla_\infty$-flat section
$\rho$ of $\W$. Unfortunately, $\rho$ is defined on a universal cover
$\tilde\M$ of $\M$, thus it gives rise to a Frobenius manifold structure on
the open subset of $\tilde\M$, where $\Phi(\rho,\cdot)$ is an isomorphism.
One easily checks that $\circ$ descends to $\M$. However,
the metric and $\tilde\nabla$ do not descend to $\M$.
In~\S\ref{primchoice} we shall construct primitive sections
on open subsets of $\M$ defined up to a scalar, then the flat structures
also descend to these open subsets of $\M$.

\subsection{Isomonodromic deformations}\label{isomonodromy}
We are going to identify $\nabla_\infty$ with some isomonodromic deformation.
Thus we shall need some generalities on isomonodromy. For more details we
refer to~\cite{Babelon}. Note that there is a standard relation between
isomonodromy and Frobenius manifolds, see~\cite{Sabbah}, but it will not be
used in our paper.

Let $\Y\to\M$ be a holomorphic fiber bundle, equipped with a family of
meromorphic flat connections on fibers (i.e., a relative connection). This
family is said to be \emph{isomonodromic} if it can be extended to a flat
absolute meromorphic connection on $\Y$. It is known that the family of connections with
regular singularities is isomonodromic if and only if the monodromy does not
change.

\subsubsection{Universal isomonodromy for line bundles}\label{universal}
We shall need a baby version of isomonodromy, namely, isomonodromy for
\emph{line bundles\/} with singular connections. In this case a family of meromorphic
connections is isomonodromic if and only if the monodromy does not change
(even if connections have irregular singularities).

Let $\M_{g,n}$ be the moduli space of triples $(X,p,x)$, where
$X$ is a curve (which we assume smooth complete over $\C$),
$p\in X$, $x$ is a coordinate to order $n$ at $p$ (that is, an $n$-jet of a coordinate). Let
$Conn_n\to\M_{g,n}$ be the moduli space of line bundles with
connections with a pole of order at most $n$ at $p$ and no
other poles. One defines the \emph{universal isomonodromic
connection\/} on this fibration by the following requirements:
a family is isomonodromic if (1) the monodromy representation is
constant and (2) the $x$-expansion of the polar part of the connection
at $p$ does not change. The existence and uniqueness of such a
connection follows easily from the Riemann--Hilbert correspondence.

\subsection{Identifying $\nabla_\infty$}
Let $m=(X,p,f_0)\in\A$. Let $\bar\A$ be the formal completion
of $\A$ at $m$. Recall that in~\S\ref{universal} we defined a
fibration $Conn_n\to\M_{g,n}$. Let $Conn_n^0\to\M_{g,n}$ be its
part corresponding to degree zero line bundles and let $\Conn
Conn^0_n\to\M_{g,n}$ be its relative universal cover. In other
words,
\begin{equation*}
    \Conn Conn^0_n=\bigsqcup_{(X,p,x)}\Conn Conn^0_{X,p,n},
\end{equation*}
where $\Conn Conn^0_{X,p,n}$ is defined in~\S\ref{32}.
Consider the diagram
\begin{equation}\label{pullback}
\begin{CD}
\phi_*\D|_{\bar\A\times\{z=\infty\}}@>>> \Conn Conn^0_n\\
@VVV @VVV\\
\bar\A @>>> \M_{g,n}
\end{CD}
\end{equation}
The lower map is given in the following way: $f$ in the neighborhood of $p$
gives a polar part of order $n$. This gives a coordinate to order $n$. In
other words, the coordinate is $x=f^{-\frac1n}$. Precisely, this coordinate
is defined up to a multiplication by a root of unity. We fix such a choice
(this is why we restrict this map to $\bar\A$).
It follows from Lemma~\ref{isomorhism} that this is a pull-back
diagram. It follows from \S\ref{universal} that the right
fibration in~(\ref{pullback}) is equipped with isomonodromic
connection.
\begin{lemma}\label{lm:nablainfty}
The upper map respects connections. Thus $\nabla_\infty$ is the pull-back of
the isomonodromic connection.
\end{lemma}
We delegate the proof to the Appendix.

We have calculated $\nabla_\infty$ on $\phi_*\D$. Now we want
to understand $\nabla_\infty$ on $\V$. To avoid confusion, we shall denote the connection
on $\V$ by $\nabla^\V$ until the end of this subsection. Clearly,
$\nabla_z=\nabla^\V_z$ for $z\ne\infty$.

\begin{lemma}\label{SplitIsoConnection}
{\rm(}a{\rm)} $\nabla^\V_\infty$ respects splitting~(\ref{splitting}).\\
{\rm(}b{\rm)} The restriction of $\nabla^\V_\infty$ to the universal cover of
the space of non-singular connections with trivial $a$-monodromy is the cover
of the isomonodromic connection.\\
{\rm(}c{\rm)} The restriction of $\nabla^\V_\infty$ to the space of trivial
bundles is described as follows: a family $(\Oo,d_X+\rho)$ is flat if and
only if the $a$-periods of $\rho$ are constant and the $x$-expansion of the
polar part of $\rho$ is constant.
\end{lemma}
\begin{proof}
(b) Follows from the previous lemma because we do not modify $\phi_*\D$ along
this subspace.\\
(a) By~(b) it is enough to prove that the space of trivial bundles is
preserved. It follows from the fact that the sections of $\V$ corresponding
to trivial bundles at
$z=\infty$ correspond to the sections of $\phi_*\D$ that vanish at $z=\infty$.\\
(c) Let $\mathcal C=(\Oo,d_X+\rho)$ be a family with constant $a$-monodromy
and a constant $x$-expansion of the polar part. Let $\mathcal C'$ be a
section of $\phi_*\D|_{z=\infty}$ covering a family of non-singular
connections such that its $a$-monodromy is trivial and the $b$-monodromy is
reciprocal to that of $\mathcal C$. Then the projection of $\mathcal
C''=\mathcal C+\mathcal C'$ to $Conn_n$ is isomonodromic in the sense
of~\S\ref{universal}. Hence it can be extended
to a $\nabla$-flat section $\tilde{\mathcal C}''$ of $\phi_*\D$ in the
neighborhood of $z=\infty$. Then $\tilde{\mathcal C}''/z$ is also
$\nabla$-flat. It can be viewed as a section of $\V$, equal
to~$(\Oo,d_X+\rho)$ at $z=\infty$. It remains to recall that $\nabla^\V_\infty$
is defined as $\lim_{z\to\infty}\nabla^\V_z$.
\end{proof}

\subsection{Deformation theory}\label{multiplication}
Denote the complex
\[
\begin{CD}
0\to\T(-1)\xrightarrow{\times\omega}\Oo(n)\to0 \\
\end{CD}
\]
by $\K$ and take $m=(X,p,f_0)\in\A$. Then the restriction $\K|_m$ governs
deformations of $m$, where only deformations of $f_0$ that preserve periods
are allowed. Thus $\K|_m$ is the deformation complex of $\A$. We give more details in~\S\ref{proofresidue}.

\begin{remark}
$\K$ is actually a dg-Lie algebra. This gives a relation of our work
to~\cite{KonBar}, where a similar complex is used (in a higher-dimensional
situation). See also~\cite{Merkulov}.
\end{remark}

Consider the residue of $\nabla$ at $z=0$, denote it by $\Phi$.
We want to give a cohomological interpretation of $\Phi$.
Consider a morphism of complexes
\begin{equation}\label{CD}
\begin{CD}
\T(-2) @>\times\omega\oplus\times\omega>> \Oo(n-1)\oplus\Oo(n-1)
@>(\times\omega)-(\times\omega)>> \Omega(2n)\\
@VVV @V-VV @VV=V \\
0@>>>\Oo(n-1) @>\times\omega>>\Omega(2n)
\end{CD}
\end{equation}
The upper complex is $\K\otimes(\Oo(-1)\xrightarrow{\times\omega}\Omega(n))$,
while the lower complex is naturally quasi-isomorphic to
$\Oo(-1)\xrightarrow{\times\omega}\Omega(n)$ shifted by one. Thus we get a map
\begin{equation*}
    \T\hat\A\otimes R^1\phi_*(\Oo(-1)\xrightarrow{\times\omega}\Omega(n))\to
    R^1\phi_*(\Oo(-1)\xrightarrow{\times\omega}\Omega(n)).
\end{equation*}
\begin{lemma}\label{residue}
This map coincides with $\Phi$.
\end{lemma}
Again we delegate the proof to the Appendix.

Notice that $\K$ is exact except at zeros of $\omega=d_Xf_0$,
therefore for a fixed $m\in\hat\A$ the complex is naturally
quasi-isomorphic to $\oplus_s\Oo_{q_s}$, where $q_s$ are zeros
of $\omega$ (if $\omega$ has multiple zeros, they should be
viewed as schemes with nilpotents). Similarly, $\D|_{z=0}$ is
quasi-isomorphic to $\oplus_s(\Omega_X)_{q_s}$.

It follows easily from the lemma above that under this
identification $\Phi$ becomes the componentwise multiplication.
Thus $\Phi(\rho_m,\cdot)$ is an isomorphism for a generic $\rho_m\in\V_m$. It follows that taking different primitive sections
we can construct twisted Frobenius structures on open subsets
covering $\hat\A$ (see~\S\ref{WDVV}).

The tangent sheaf multiplication also has a cohomological interpretation,
namely, there is a natural map $\K\otimes\K\to\K\otimes\Oo(n)$, similar
to~(\ref{CD}).
\begin{lemma}
The induced map on cohomology coincides with $\circ$. In particular $\circ$
does not depend on the choice of a primitive section.
\end{lemma}
\begin{proof}
Follows from associativity of cohomological multiplication.
\end{proof}
One would like to have a canonical choice of a primitive section; this will
be discussed in~\S\ref{PrimitiveSections}.

\subsection{Metric}
We shall sketch a construction of a bilinear product on $\V$. We start
with the natural $\phi^\cdot\Oo_{\hat\A\times\Pro^1}$-linear map of complexes on $\X\times\Pro^1$
\begin{equation*}
\D(1-n)\otimes
a^*\D\to(\Oo(-n-1)\xrightarrow{d}\Omega)\boxtimes\Oo_{\Pro^1}(1).
\end{equation*}
For $z\ne\infty$ it gives rise to a $\nabla$-flat map:
$\langle\cdot,\cdot\rangle:\phi_*\D\otimes a^*\phi_*\D\to\Oo_{\Pro^1}(1)$ (it
is a priori singular at $z=\infty$ because $\D(1-n)$ is not quasi-isomorphic to $\D$
at $z=\infty$). This form is \emph{$a$-skew-symmetric}.

We claim that the restriction of this metric to $\V$ vanishes at $z=\infty$.
Indeed, let~$s_1$, $s_2$ be sections of $\V$. Then, by definition of $\V$,
$s_i|_{z=\infty}$ are tangent vectors to the space of non-singular connection
with zero $a$-periods. One can check that our form agrees with the natural
symplectic form on the space of connections. Notice that the symplectic form
on non-singular connections can be identified with the intersection form on
$H^1(X,\C)$. It remains to recall that $\Delta$ is isotropic.

Thus composing $\langle\cdot,\cdot\rangle|_\V$ with multiplication by $z$ we
get an \emph{$a$-symmetric} bilinear form on $\V$. We leave it to the reader
to check that the form is non-degenerate.

\begin{remark}
We have to multiply by $z$ because we identified $\Oo_{\Pro^1}(1)$ with
functions with a pole at $z=0$ in~(\ref{global}).
\end{remark}

\subsection{Relation with the construction of Krichever}
The above Frobenius structure is equivalent to that of~\cite{Krichever1,
Krichever2}. To give a bridge between these papers and ours we present here a
different point of view on $\V$.

Fix $m=(X,p,f_0,\Delta)\in\hat\A$ and fix some set of disjoint closed curves
$a_i$ ($i=1,\ldots,g$) on $X$ representing $a$-cycles. Assume that
$p\notin\cup a_i$. Let $\chi$ be a 1-form on $X$ with the only pole of order
at most $n$ at $p$ and \emph{jumps\/} $\lambda_i\omega$ along $a_i$, where
$\lambda_i\in\C$ ($\omega=d_Xf$ as before).

\begin{lemma}
For every $z\in\Pro^1$ there is a natural isomorphism between the set of such
forms with jumps and $\V|_{m,z}$.
\end{lemma}
\begin{proof}
Let us first consider $z\ne\infty$. The idea is that a form with jumps can
be, in some sense, viewed as a \v{C}ech cocycle of $\D$.

Precisely, we take a cover $\U_j$ of $X$ such that there is a
unique $j$ with $p\in\U_j$. Let $(\alpha_j,s_{jk})$ be a
cocycle (cf.~\S\ref{calc}), representing a class in $\Hyper^1(X,\D)$ so that
$\alpha_j\in\Gamma(\U_j,\Omega(n))$,
$s_{jk}\in\Gamma(\U_j\cap\U_k,\Oo)$. The cocycle conditions are
\[
\alpha_j-\alpha_k=d_Xs_{jk}+\frac\omega z s_{jk},\qquad
s_{jk}+s_{kl}+s_{lj}=0,\qquad s_{jk}=-s_{kj}.
\]
Hence $s_{jk}$ represents a class in $H^1(X,\Oo)$. It is easy
to check that for every such a cocycle $s_{jk}$ we can find
holomorphic functions $h_j$ \emph{with constant jumps along
$a_i$'s\/} such that $s_{jk}=h_j-h_k$. This functions are
unique if we require $h_j(p)=0$, if $p\in\U_j$. It follows that
\[
\alpha_j-d_Xh_j-\frac\omega zh_j
\]
patch together to a 1-form $\alpha$ with jumps along $a_i$. We leave it to the
reader to check that this map from the cohomology group to forms with jumps
is an isomorphism.

For $z=\infty$ one gives an isomorphism by the following conditions\\
\indent\noindstep If $\lambda_i=0$ for all $i$, then $\chi$ corresponds to
$(\Oo,d+\chi)$.\\
\indent\noindstep If $\chi$ has no pole at $p$, then consider the local
system given by the transition functions $\exp(\lambda_i)$ (in particular its
$a$-monodromy is trivial). This map to the space of local systems on $X$ is
then lifted to a map to $\V|_{z=\infty}$.
\end{proof}

The reader should compare this description with
Lemma~\ref{SplitIsoConnection}. It follows that in terms of forms with jumps
$\nabla_\infty$ can be described very easily: a family is flat if $a$-periods
do not change, $\lambda_i$'s do not change and the $x$-expansion of $\chi$ at
$p$ does not change.

Krichever fixes a basis in $H_1(X,\Z)$ and a choice of $f^{-1/n}$. It allows
him to construct a flat basis in the space of forms with jumps (denoted by
$\Omega_a$). See~\S7.5 and~\S7.2 of~\cite{Krichever1}.

\section{Multi-point generalizations}\label{sectmulti}
In this section we generalize the above constructions to the case of abelian integrals with many poles. To define such an integral we need an extra ``reference point''.

Let $X$ be as before, $p_0,\ldots,p_k$ be distinct points of $X$. Let $n_1,\ldots,n_k$ be non-negative integers and consider the divisors $\divisor=p_0+\sum_{i=1}^k n_ip_i$, $\divisor'=\sum_{i=1}^k p_i$; we assume that $\sum n_i+k\ge2$, $k\ge1$. By an \emph{abelian integral\/} with polar divisor $\divisor$ we mean a function $f$ on the universal cover of $(X,p_0)$ such that $df$ descends to a meromorphic differential on $X$ with a pole of order exactly $n_i+1$ at $p_i$ for $1\le i\le k$ and no other poles (note that there is no pole at $p_0$). Let $\A_{g,n_1,\ldots,n_k}$ be the corresponding moduli space of collections $(X,p_0,\ldots,p_k,f)$.

For fixed $(X,p_1,\ldots,p_k)$ we get a period map $H_1(X\setminus\{p_1,\ldots,p_k\},\Z)\to\C$. As before, identifying the homology groups for nearby curves, we get a period foliation on
$\A_{g,n_1,\ldots,n_k}$. Fix one of the leaves $\A$.

\begin{remark}
Note that $\A_{g,n_1,\ldots,n_k}=\A'\times\C$,  where $\A'$ is the corresponding moduli space of abelian differentials. The splitting is given by the abelian integrals, which vanish at $p_0$.
\end{remark}

For $(X,p_0,\ldots,p_k,f)\in\A$ we have a natural map $H_1(X\setminus\{p_1,\ldots,p_k\},\Z)\to H_1(X,\Z)$ and we pull the intersection form back via this map to a degenerate alternating form on $H_1(X\setminus\{p_1,\ldots,p_k\},\Z)$.
Let $\Delta\subset H_1(X\setminus\{p_1,\ldots,p_k\},\Z)$ be maximal isotropic with respect to this form. Then $\Delta$ has rank $g+k-1$.
Denote by $\hat\A$ the module space of collections $(X,p_0,\ldots,p_k,f,\Delta)$, where $(X,p_0,\ldots,p_k,f)\in\A$, $\Delta$ is as above.

Again we denote by $\phi:\X\to\A$ the natural projection from the universal curve, by $\tilde p_i$ the universal divisors on $\X$ corresponding to $p_i$. Set $\tilde\divisor=\tilde p_0+\sum_1^k n_i\tilde p_i$, $\tilde\divisor'=\sum_1^k\tilde p_i$. Analogously to~\eqref{global}, we set $\omega=d_Xf$ and consider the complex
\[
    \D:=\left(\Oo(-\tilde\divisor')\boxtimes\Oo_{\Pro^1}\xrightarrow{d_X+\frac\omega z}
    \Omega(\tilde\divisor)\boxtimes\Oo_{\Pro^1}(1)\right).
\]
Similarly to Lemma~\ref{structure}, we prove that the restriction of $\phi_*\D$ to a fixed point of $\A$ is isomorphic to the vector bundle
$\left(\oplus_1^{g+k-1}\Oo_{\Pro^1}\right)\oplus
\left(\oplus_1^{g+\sum n_i}\Oo_{\Pro^1}(1)\right)$.

Similarly to~\S\ref{universal} we introduce the moduli space $\M_{g;n_1,\ldots,n_k}$ parameterizing
collections $(X,p_0,\ldots,p_k,x_1,\ldots,x_k)$, where
$X,p_0,\ldots,p_k$ are as before, $x_i$ is a coordinate to order $n_i$ at $p_i$.
Note that if $n_i=0$, then $x_i$ is just the point $p_i$ itself.

Let $Conn^0_{n_1,\ldots,n_k}\to\M_{g;n_1,\ldots,n_k}$ be the moduli space of degree zero line bundles
trivialized over $\divisor'$ with
connections with a pole bounded by $\divisor$. Let $\Conn
Conn^0_{n_1,\ldots,n_k}\to\M_{g;n_1,\ldots,n_k}$ be the relative universal cover.

Let $m\in\A$. Let $\bar\A$ be the formal completion
of $\A$ at $m$. We have
\begin{lemma}
Let $\omega$ be a meromorphic form on a disc with a pole of order $n\ge2$ at the center. Then there is a unique up to multiplication by a root of unity $n-1$-jet of coordinate $x$ such that $\omega=(x^{-n}+ax^{-1})\,dx+O(1)$, where $O(1)$ is a regular form.
\end{lemma}
\begin{proof}
Let $\omega_0$ be a form with a simple pole and the same residue as $\omega$. Then
\[
    x=c\left(\int(\omega-\omega_0)\right)^{\frac{-1}{n-1}},
\]
where $c$ is a suitable constant.
\end{proof}
Applying this Lemma at every $p_i$ with $n_i>0$, we get a map $\psi:\bar\A\to\M_{g;n_1,\ldots,n_k}$.

Similarly to~\eqref{pullback} we have a pull-back diagram
\begin{equation*}
\begin{CD}
\phi_*\D|_{\bar\A\times\{z=\infty\}}@>>> \Conn Conn^0_{n_1,\ldots,n_k}\\
@VVV @VVV\\
\bar\A @>\psi>> \M_{g;n_1,\ldots,n_k}
\end{CD}
\end{equation*}
Again the right fibration is equipped with isomonodromic
connection and the upper map respects connections. Thus $\nabla_\infty$ is the pull-back of
the isomonodromic connection. As before we get a splitting
\[
\phi_*\D|_{\bar\A\times\{z=\infty\}}=\widetilde{Bun}\oplus Conn_{triv},
\]
where $\widetilde{Bun}$ is the relative universal cover of the moduli space of bundles with trivialization on $\divisor'$ and non-singular connections with zero a-periods, and $Conn_{triv}$ is the moduli space of connections on the trivial bundle with standard trivialization on $\divisor'$. Then we define $\V$ as the vector bundle on $\hat\A$ whose sheaf of sections is
\[
\{s\in\phi_*\D: s|_{z=\infty}\in\widetilde{Bun}\}.
\]
We have a multi-pole version of Theorem~\ref{identifications}: $\V=p_1^*\W$ for some
vector bundle $\W$ on $\hat\A$ and
$(\W,\nabla)$ is a pencil of connections on $\hat\A$ in the sense of
Definition~\ref{pencil}. Finally, the multiplication and the metric have a cohomological description similar to the single-point case.

\section{Primitive sections}\label{PrimitiveSections}
\subsection{A general strategy} Let $g\ge0$ and $2\le k\le n$ be integers. We start with a complete smooth curve $Y$ of genus $g$, a point $q\in Y$ and differentials $\omega_0\in H^0(Y,\Omega_Y((n+1)q))$ and $\rho_0\in H^0(Y,\Omega_Y(kq))$ having no common zeros. Then choose an abelian integral $f_0$ of $\omega_0$, and consider the local coordinate $x=f_0^{-1/n}$ near $q$.

Now consider the moduli space $\A'$, parameterizing collections $(X,p,f,\Delta,\rho)$, where $X\ni p$ is a complete curve of genus $g$, $f$ is an abelian integral such that $df\in H^0(X,\Omega_X((n+1)p))$, $f$ has the same periods as $f_0$ (more precisely, the collection is in the same leaf of the period foliation, cf.~\S\ref{mainobj}), $\Delta$ is maximal isotropic in $H_1(X,\Z)$, $\rho\in H^0(X,\Omega_X(kp))$ such that $\rho$ and $\rho_0$ have the same $a$-periods, and the coefficients of the $f^{-1/n}$-expansion of $\rho$ at $p$ are equal to the coefficients of the $x$-expansion of $\rho_0$ at $q$. By Lemma~\ref{SplitIsoConnection}(c), $\rho$ is $\nabla_\infty$-flat. Let $\A''$ be the non-empty open subset of $\A'$, where $\rho$ and $\omega$ have no common zeros. Then $\rho$ is a primitive section on $\A''$, and $\A''$ is a Frobenius manifold. We leave the multipoint generalization to the reader.

\subsection{Another choice of a primitive section}\label{primchoice}
The following choice of a primitive section is often used in the theory of integrable systems. Let $\hat\A$ be as in~\S\ref{mainobj}. We use notation of~\S\ref{notation}. Let $k$ be an integer such that $2\le k\le n$. Consider a section $\rho_k$ of
$\Omega(k)$ with the following properties\\
\indent\noindstep Its polar part at $\tilde p$ is of the form $x^{-k}dx$,
where $x=f^{-1/n}$.\\
\indent\noindstep The $a$-periods of $\rho_k$ are zero.

This form is defined locally over $\hat\A$ up to a multiplication by an $n$th root of
unity. It follows from Lemma~\ref{SplitIsoConnection}(c) that $(\Oo,d_X+\rho)$
is a flat section of $\V|_{z=\infty}$. This is a primitive section on the
open set $\hat\A_k\subset\hat\A$, where $\rho_k$ and $\omega$ have no common
zeros. Thus it gives rise to a connection on $\hat\A_k$, which depends only
on $k$. The metric is defined up to multiplication by a root of unity.

It is curious that if $n$ is even, then $\rho_{n/2}$ is defined up to sign,
so the corresponding metric is also defined canonically.

Unfortunately, we do not know whether $\hat\A_k$ is always
non-empty. We have the following partial result
\begin{lemma}
{\rm(}a{\rm)} For a generic leaf $\A\subset\A_{g,n}$ all the sets $\hat\A_k$ are non-empty.\\
{\rm(}b{\rm)} If $n>k+2g$, then $\hat\A_k$ is non-empty.
\end{lemma}
\begin{proof}
(a) It is enough to show that for a generic $(X,p,f)\in\A_{g,n}$ the corresponding $\rho$ and $df$ have no common zeros. Consider first the case $g>0$. Changing $f$ to $f+\int\alpha$, where $\alpha$ is any holomorphic differential on $X$ does not change the $k$-jet of the coordinate $x=f^{-1/n}$. Thus it does not change $\rho$, so it is enough to show that there is $\alpha\in H^0(X,\Omega_X)$ such that $df+\alpha$ and $\rho$ have no common zeros. Let $\{q_i\}$ be the set of zeros of $\rho$, let $H_i$ be the set of holomorphic differentials $\alpha$ such that $(df+\alpha)(q_i)=0$. Since $H_i$ is an affine subspace of $H^0(X,\Omega_X)$, it is enough to show that $H_i\ne H^0(X,\Omega_X)$ (because finite number of proper affine subspaces cannot cover a vector space). Thus we have to show that the natural map $H^0(X,\Omega_X(-q_i))\to H^0(X,\Omega_X)$ is not an isomorphism. Using a standard exact sequence we see, that it is equivalent to $H^1(X,\Omega_X(-q_i))=\C$, which, by Serre duality, is equivalent to $H^0(X,\Oo_X(q_i))=\C$. The last equality follows from $g>0$.

Now assume that $g=0$, so that $X=\Pro^1$, and we can assume $p=\infty$. Since $\A_{0,n}$ is irreducible, it is enough to find a single $f$ such that $df$ and $\rho$ have no common zeros. We can take $f=z^n+z$, where $z$ is the standard coordinate. Then an easy calculation shows that $\rho=-z^{k-2}\,dz$.

(b) Let $\A$ be any leaf, $(X,p,f)\in\A$. If we add a section of $H_0(X,\Oo_X((n-k)p))$ to $f$, the $k$-jet of $f^{-1/n}$ does not change, so $\rho$ does not change, hence, as before, it is enough to show that there is $g\in H_0(X,\Oo_X((n-k)p))$ such that $df+dg$ and $\rho$ have no common zeros. The rest of the proof is as above. In the last step we use that a sheaf of degree at least $2g-1$ has zero first cohomology group.
\end{proof}

\section{Appendix: Some proofs}
\subsection{Calculating $\phi_*\D$}\label{calc}
We start with the \v{C}ech calculation of $\phi_*\D$. Pick $m=(X,p,f)\in\A$, let $\bar\A$ be the completion of $\A$ at $m$. Consider
the following cover of $X$: $X=\dot X\cup{}D$, where $\dot
X=X\setminus p$, $D$ is the formal neighborhood of $p$. Let
$\dot D=\dot X\cap D$ be the punctured formal neighborhood. Let $\bar\X$ be the restriction of $\X$
to $\bar\A$. Since affine schemes have no infinitesimal
deformations, $\bar\X$ can be covered by $\dot X\times\bar\A$
and $D\times\bar\A$. To calculate $\phi_*\D$ we shall use the
\v{C}ech resolution of $\D$
\begin{equation*}
\begin{CD}
0 @>>>\Oo(-1) @>d_X+\frac\omega z>>
\Omega(n) @>>>0 @>>>0\\
@VVV @VVV @VVV @VVV @VVV\\
0 @>>>C^0(\Oo(-1))
@>>>\genfrac{}{}{0pt}{}{C^1(\Oo(-1))\oplus}{C^0(\Omega(n))} @>>>
C^1(\Omega(n)) @>>> 0
\end{CD}
\end{equation*}
Consider a section of $R^1\phi_*\D$ over $\bar\A$. It is
represented by a family of 1-cocycles
$(\alpha_1,\alpha_2,s)$.
Precisely, $\alpha_1$ is a relative 1-form on $\dot
X\times\bar\A\to\bar\A$, $\alpha_2$ is a relative
1-form on $D\times\bar\A\to\bar\A$ with a pole of order at most
$n$ along $\{p\}\times\bar\A$, $s\in\Oo_{\dot
D\times\bar\A}$. The cocycle condition is
\begin{equation}\label{cocycle}
    \alpha_1-\alpha_2=d_Xs+\frac{\omega s}z.
\end{equation}
We extend $d_X+\frac\omega z$ to an absolute connection
$\nabla=d+\frac{df}z$ as in~\S\ref{MainConstruction}.

We want to describe the Gauss--Manin connection on the first
hypercohomology sheaf directly.
First, we define a map
\begin{equation*}
    \Psi_1:\wedge^2\Omega_{\dot X\times\bar\A}\to
    \phi^*\Omega_{\bar\A}\otimes\Omega_{\dot X\times\bar\A/\bar\A}.
\end{equation*}
It is defined as follows: there is a natural surjective map
$\phi^*\Omega_{\bar\A}\otimes\Omega_{\dot
X\times\bar\A}\to\wedge^2\Omega_{\dot X\times\bar\A}$ (recall
that $\dim X=1$). Thus we can lift a section of
$\wedge^2\Omega_{\dot X\times\bar\A}$ to
$\phi^*\Omega_{\bar\A}\otimes\Omega_{\dot X\times\bar\A}$ and
then project it to $\phi^*\Omega_{\bar\A}\otimes\Omega_{\dot
X\times\bar\A/\bar\A}$. The lift is defined up to an element of
the second symmetric power of $\phi^*\Omega_{\bar\A}$, thus the
projection does not depend on the lift. Similarly, we can
define maps
\begin{equation*}
\begin{split}
    \Psi_2:\wedge^2\Omega_{D\times\bar\A}\to
    \phi^*\Omega_{\bar\A}\otimes\Omega_{D\times\bar\A/\bar\A}.\\
    \Psi_{12}:\wedge^2\Omega_{\dot D\times\bar\A}\to
    \phi^*\Omega_{\bar\A}\otimes\Omega_{\dot D\times\bar\A/\bar\A}.
\end{split}
\end{equation*}
Now let us lift $\alpha_i$ to an absolute 1-form
$\tilde\alpha_i$ such that $\tilde\alpha_2$ has a pole of order
at most~$n$ on $\{p\}\times\bar\A$.
It follows from~(\ref{cocycle}) that
\begin{equation}\label{beta}
    \tilde\alpha_1-\tilde\alpha_2=\nabla s+\beta,
\end{equation}
where $\beta$ is a section of $\phi^*\Omega_{\bar\A}$. The
Gauss--Manin connection is given by
\begin{equation}\label{GaussManin}
    \nabla(\alpha_1,\alpha_2,s)=(\Psi_1(\nabla\tilde\alpha_1),
    \Psi_2(\nabla\tilde\alpha_2),\beta).
\end{equation}
\begin{lemma}
\noindstep\label{One} The 1-cochain~(\ref{GaussManin}) is a
cocycle for
$\D\otimes\Oo(n+1)\otimes\phi^*\Omega_{\bar\A}$.\\
\noindstep\label{Two} If we change the lifts
$\alpha_i\rightsquigarrow\tilde\alpha_i$,
then~(\ref{GaussManin}) changes by a coboundary.\\
\noindstep\label{Three} $\nabla$ satisfies the Leibnitz rule
with respect to multiplication of cocycles by functions on
$\bar\A$.
\end{lemma}
\begin{proof}
We have a well-defined relative differential
$d_X:\phi^*\Omega_{\bar\A}\to
\phi^*\Omega_{\bar\A}\otimes\Omega_{X\times\bar\A/\bar\A}$. The
proof of~(\ref{One}) and~(\ref{Two}) is based on the following:
\emph{if $\gamma$ is a section of $\phi^*\Omega_{\bar\A}$, then
$\Psi_i(\nabla\gamma)=d_X\gamma+\frac1z\gamma\otimes\omega$.}
The proof of this statement is left to the reader, let us prove
(\ref{One}). It follows from~(\ref{beta}) that
$\nabla\tilde\alpha_1-\nabla\tilde\alpha_2=\nabla\beta$. Now
let us apply $\Psi_{12}$ to both sides, it gives
$\nabla\alpha_1-\nabla\alpha_2=d_X\beta+\frac1z\beta\otimes\omega$
and the cocycle condition follows.

We leave the proofs of~(\ref{Two}) and~(\ref{Three}) to the
reader.
\end{proof}
It follows that $\nabla$ descends to a connection on
$\phi_*\D$. A little difficulty is that we get a cocycle of
$\D\otimes\Oo(n+1)\otimes\phi^*\Omega_{\bar\A}$ instead of that of $\D\otimes\phi^*\Omega_{\bar\A}$. Fortunately,
there is a natural quasi-isomorphism (compare
with~\S\ref{Dmodules})
\begin{equation}\label{quasiisomorphism}
    \D\hookrightarrow\D\otimes\Oo(n+1)
\end{equation}
for $z\ne\infty$. For $z=\infty$ one needs to extend the
cocycle to a neighborhood of $z=\infty$ first, apply $\nabla$,
and then pass to the limit as $z\to\infty$ to evaluate
$\nabla_\infty$. This will be done in the next subsection.

\subsection{Proof of Lemma~\ref{lm:nablainfty}}
Consider an isomonodromic family of connections on $\bar\A$
given by a family of cocycles
$(\alpha_1,\alpha_2,\exp(s))$ as
in \S\ref{calc} (now we have $z=\infty$). Here
$\exp(s)$ is a cocycle defining the line bundle. It
has a single-valued logarithm because the bundle has degree
zero. We need to show that this family is $\nabla_\infty$-flat.
The cocycle condition~(\ref{cocycle}) becomes
\begin{equation*}
    \alpha_1-\alpha_2=d_Xs.
\end{equation*}
Since the family is isomonodromic, the forms $\alpha_i$ can be
extended to absolute \emph{closed\/} forms $\tilde\alpha_i$.
The condition that the $x$-expansion of the polar part does not
change in the family can be written in the following way
\begin{equation}\label{Alpha}
    \tilde\alpha_2=h(x^{-1})dx+\gamma(\epsilon,x),
\end{equation}
where $\epsilon$ is a coordinate on $\bar\A$, $h$ is a polynomial with constant coefficients, $\gamma$
is a 1-form regular on $\{p\}\times\bar\A$. A~simple calculation in
local coordinates shows that, changing the lift
$\alpha_2\rightsquigarrow\tilde\alpha_2$ if necessary, we can
assume that $\gamma/x$ has a logarithmic pole on $\{p\}\times\bar\A$.

Since $f=x^{-n}$, we have $h(x^{-1})dx\wedge df=0$. Now we want
to extend the cocycle to the neighborhood of $z=\infty$. To
this end we write $\omega s=\sigma_1-\sigma_2$, where
$\sigma_1$ is a relative 1-form on $\dot X\times\bar\A$,
$\sigma_2$ is a relative 1-form on $D\times\bar\A$ with at most
simple pole on $\{p\}\times\bar\A$. This is always possible, since
$H^1(X,\Omega_X(p))=0$. It is easy to see that
$(\alpha_1+\frac{\sigma_1}z,\alpha_2+\frac{\sigma_2}z,s)$
satisfies the cocycle condition~(\ref{cocycle}).

We extend $\sigma_i$ to an absolute form $\tilde\sigma_i$. One checks that we can choose $\tilde\sigma_2$ so that it has at most simple pole on $\{p\}\times\bar\A$ and $\tilde\sigma_2\wedge df=0$.
Using~(\ref{Alpha}) we get
\begin{equation*}
\nabla\left(\tilde\alpha_2+\frac{\tilde\sigma_2}z\right)=
\left(d+\frac{df}z\right)\left(\tilde\alpha_2+\frac{\tilde\sigma_2}z\right)=
\frac{df}z\wedge\gamma+\frac{d\tilde\sigma_2}z=O\left(\frac1z\right).
\end{equation*}
A similar argument, shows that
$\nabla\left(\tilde\alpha_1+\frac{\tilde\sigma_1}z\right)=O(1/z)$.
Applying $\Psi$ we see that
\begin{equation}\label{notapriori}
    \nabla\left(\alpha_1+\frac{\sigma_1}z,\alpha_2+\frac{\sigma_2}z,s\right)=(0,0,\beta)+O\left(\frac1z\right),
\end{equation}
for some $\beta$. Note also that since $\gamma/x$ has
logarithmic pole on $\{p\}\times\bar\A$, (\ref{notapriori}) is a cocycle
of $\D$ (a priori it is a cocycle of $\D(n)$). Therefore we do
not need to invert the
quasi-isomorphism~(\ref{quasiisomorphism}) (which could have
altered the behavior at $z=\infty$). Taking limit as
$z\to\infty$, we get a cocycle of the form $(0,0,\beta)$.
The cocycle condition shows that $d_X\beta=0$, and we see
that this cocycle is a coboundary, so the family $(\alpha_1,\alpha_2,s)$ is a
flat section of $\nabla_\infty$.

\subsection{Proof of Lemma~\ref{residue}}\label{proofresidue}
Let us fix $m=(X,p,f_0,\Delta)\in\hat\A$ and give an explicit
description of the map $\Hyper^1(\K|_m)\to\T_m\hat\A.$ Consider a
cocycle $(h_1,h_2,\tau)$ of $\K|_m$, where $h_1\in\Oo_{\dot X}$,
$h_2\in\Oo_D(n)$, $\tau\in\T(\dot D)$. Then $\tau$ represents a
class of $H^1(X,\T)$, thus it gives rise to an infinitesimal
family of curves $\hat X\to\Spec\C[\epsilon]/\epsilon^2$. A function
$\hat f$ on $\hat X$ is given by the conditions: $\hat f|_{\dot
X\times\Spec\C[\epsilon]/\epsilon^2}=f_0+\epsilon h_1$, $\hat
f|_{D\times\Spec\C[\epsilon]/\epsilon^2}=f_0+\epsilon h_2$. The
periods of $\hat f$ are constant in the family because $h_1$ is
a single-valued function on $\dot X$. Thus we have assigned a
vector of $\T_m\hat\A$ to $(h_1,h_2,\tau)$. We leave it to the
reader to check that it indeed gives an isomorphism of
$R^1\phi_*\K$ with the tangent sheaf of $\hat\A$.

Let us take a \v{C}ech cocycle $(\alpha_1,\alpha_2,s)$ of
$\D|_{m,z=0}$. We have $\alpha_1\in\Omega_{\dot X}$,
$\alpha_2\in\Omega_D(n)$, $s\in\Oo_{\dot D}$. Let us also take
a cocycle $(h_1,h_2,\tau)$ of $\K|_m$, where $h_1\in\Oo_{\dot X}$,
$h_2\in\Oo_D(n)$, $\tau\in\T(\dot D)$. This cocycle represents
some $\xi\in\T_m\hat\A$. The product of cocycles is given by
$(h_1\alpha_1,h_2\alpha_2,\tau\alpha_1+h_2s)$. The cocycle
condition for $(\alpha_1,\alpha_2,s)$ is given by
$\alpha_1-\alpha_2=\omega s$, hence $s$ is uniquely determined
by $\alpha_1$ and $\alpha_2$. Again, we view $\bar\X$ as glued
from $\dot X\times\bar\A$ and $D\times\bar\A$. We extend $\alpha_1$
to a 1-form $\tilde\alpha_1$ on $\dot X\times\bar\A$ using this direct
product structure. Similarly we extend $\alpha_2$ to $\tilde\alpha_2$ and
$h_i$ to $\tilde h_i$. Let us lift $\xi$ to a vector field $\tilde\xi_1$
along $\dot X\times\{m\}$ and to $\tilde\xi_2$ along $D\times\{m\}$. In
particular we have $\langle\tilde\xi_i,\tilde\alpha_i\rangle=0$.

Unwinding the definition of $\nabla$, we see that $\Phi$ is
given on cocycles by the formula
\begin{equation*}
    \Phi(\xi,(\alpha_1,\alpha_2,s))= (\langle\xi,\Psi_1(df\wedge\tilde\alpha_1)\rangle,
    \langle\xi,\Psi_2(df\wedge\tilde\alpha_2)\rangle,?).
\end{equation*}
Since the last entry of the cocycle is uniquely determined by
the others, it is enough to prove that
\begin{equation*}
    \langle\xi,\Psi_i(df\wedge\tilde\alpha_i)\rangle=h_i\alpha_i.
\end{equation*}
Using the identification of cohomology of $\K|_m$ with the tangent
space of $\hat\A$ at $m$, one finds that
$\langle\tilde\xi_i,df\rangle=\partial_{\tilde\xi_i}(f)=h_i$.
It follows that
\begin{equation*}
\langle\xi,\Psi_i(df\wedge\tilde\alpha_i)\rangle=
\langle\tilde\xi_i,df\wedge\tilde\alpha_i\rangle=
\langle\tilde\xi_i,df\rangle\tilde\alpha_i-\langle\tilde\xi_i,\tilde\alpha_i\rangle
df=h_i\tilde\alpha_i.
\end{equation*}
The first equality is a ``consistency property'' of $\Psi_i$.

\end{document}